\documentclass[12pt]{article}
\usepackage{amssymb}
\usepackage{graphicx}
\begin{document}
\title{Incompressible representations of the Birman-Wenzl-Murakami algebra}
\author {  Vincent Pasquier.}
\date{}
\maketitle \hskip-6mm { Service de Physique Th\'eorique, C.E.A/
Saclay, 91191 Gif-sur-Yvette, France.}

\begin{abstract}
We construct a representation of the Birman-Wenzl-Murakami algebra
acting on a space of polynomials in $n$ variables vanishing when
three points coincide. These polynomials are closely related to
the Pfaffian state of the Quantum Hall Effect and to the
components the transfer matrix eigenvector of a O(n) crossing loop
model.

\end{abstract}
\maketitle

\section{Introduction and Conclusion}

This paper continues a preceding one \cite{incompressible} where
we have studied the deformations of certain wave functions of the
quantum Hall effect. There, the deformations were related to
representations of the Temperley-Lieb algebra. Here, we consider
yet another class whose deformation yields a representation of the
Birman-Wenzl-Murakami (B.M.W.) algebra \cite{birman}\cite{murak}.

Read and Rezayi \cite{Read-R1}\cite{Read-R2} have  initiated the
study quantum Hall wave functions related to parafermion currents
in two-dimensional conformal field theory. They are symmetric (or
antisymmetric) polynomials in $z_1,...,z_n$ vanishing
quadratically as any $k+1$ particles come to the same point. They
model a quantum Hall state where the particles are gathered in
clusters of $k$ and are expected to obey non-Abelian statistics.

Here, we are mainly considering the ground state which is the
minimal degree nonzero symmetric polynomial satisfying the
constraint. For a number of variables multiple of $k$, $n=km$, it
is obtained as a symmetrized sum of terms obeying the constraint
and having the factorized expression:
\begin{eqnarray}
 P(z_1,...,z_n)=\prod_{i=1}^n\prod_{p=1}^{m-1}b_{ii+p}.
 \label{read-rezayi}
\end{eqnarray}
with $b_{ij}=z_i-z_j$, and the cyclic identification of the
indices $(z_i=z_{i+n})$ is assumed.

Feigin, Jimbo, Miwa and Mukhin \cite{FJMM1} have independently
considered these wave functions. Moreover, they have introduced a
q-deformation obeying the cancellation condition \cite{FJMM2}:
\begin{eqnarray}
P(z_1,z_2,...,z_{k+1})=0\ \ {\rm for\ } \ z_{i+1}=q^{2i}z_1,\
{1\le i\le k}.
 \label{FJMM}
\end{eqnarray}
In particular, the ground state is obtained  by substituting $
b_{ij}\to qz_i-q^{-1}z_j$ in (\ref{read-rezayi}).

This paper is a step towards the generalization of the condition
(\ref{FJMM}) when the polynomials are no longer symmetrical.
Namely, we consider polynomials in an infinite set of variables
$z_i$, subject to the identification $z_{i+n}=q^{2k}p^2z_i$, and
obeying the cancellation condition:
\begin{eqnarray}
P(z_{i_1}=z,z_{i_2}=q^2z,...,z_{i_{k+1}}=q^{2k}z)=0\ \ {\rm for\ }
 {i_1<i_2...<i_{k+1}}.
 \label{FJMM-new}
\end{eqnarray}
The wave function (\ref{read-rezayi})  satisfies (\ref{FJMM-new}),
but is not the unique ground state. The claim is that the space of
the ground states can be organized into a representation of the
braid group.
Here, we concentrate on the case $k=2$. We show that the
polynomials of the minimal degree obeying (\ref{FJMM-new})  form
an irreducible representation of the Birman-Wenzl-Murakami
algebra.

\smallskip\smallskip

>From a different perspective, Di Francesco and Zinn-Justin
\cite{Pdf2} and Knutson and Zinn-Justin \cite{Zinn} have obtained
the ground state eigenvector of a $O(1)$ loop model on the
cylinder and derived from it combinatoric properties of certain
algebraic varieties, in connection with conjectures of De Gier and
Nienhuis \cite{degier}. Here, when $k=2$ and the cyclic condition
$q^{2k}=p^{-2}$ is satisfied, the degenerate ground states are the
components of the Perron-Frobenius eigenvector of a statistical
mechanics transfer matrix. Their sum is a symmetrical polynomial
given by a deformation of the Pfaffian state \cite{moore-read}:

\begin{eqnarray}
 {\rm P f}\left({1\over qz_i-q^{-1}z_j }-{1\over qz_j-q^{-1}z_i }\right)
 \prod_{i<j}^{n}{ (qz_i-q^{-1}z_j)(qz_j-q^{-1}z_i)\over z_i-z_j}.
\label{Pfaffien-pfaf}
\end{eqnarray}

\smallskip

The correspondence between the polynomials and the braid group
relies on a representation of the permutations inherited from a
solution of the Yang-Baxter equation. The starting point is a wave
function which is transformed in the same way by permuting the
coordinates or by acting with the Yang-Baxter matrix. The approach
is similar, but technically more involved than
\cite{incompressible} because the polynomials obeying
(\ref{FJMM-new}) decompose into three pieces instead of two:
\begin{eqnarray}
P=S_0+(qz_i-q^{-1}z_{i+1})(S_1+(pz_i- p^{-1}z_{i+1})S_2).
 \label{3pieces}
\end{eqnarray}
The three polynomials of the decomposition satisfy the constraint,
and the polynomials $S_i$ are symmetrical under the exchange of
$z_i$ and $z_{i+1}$. The generators of the B.M.W. algebra can be
identified with the projectors onto the three pieces of
(\ref{3pieces}).

\smallskip

The paper is divided into two sections. In the first one we recall
the definition of the B.M.W. algebra, we define the representation
which we need and we exhibit a basis in terms of tangles. In the
second one, we consider the polynomial space introduced above and
we identify a basis of these polynomials with the tangle basis.
The appendix \ref{brauer-appendice} studies the Brauer algebra
which is the rational limit  of the B.M.W. algebra. We show in the
appendix \ref{T.L.appendice} that the Temperley-Lieb structure is
intimately related to $k=2$.

We hope to return to the $k>2$ cases and to the higher degree
polynomials (excited states) elsewhere.

\bigskip

\subsection{Acknowledgements}

I have benefited from discussions with  Philippe Di Francesco,
Vincent Lafforgue and N. Read.

\section{Birman-Wenzl-Murakami algebra}

The B.M.W. algebra $\mathcal{B}_n$ is generated by braid
generators $t_i$ and Temperley Lieb generators $e_i,$ for $1\le
i\le n-1,$ obeying the relations:
\begin{eqnarray}
 &&t_it_{i+1}t_i=t_{i+1}t_it_{i+1},\ \ t_it_j=t_jt_i,\ {\rm
 if}\ |i-j|\ge 2,\cr
 &&e_ie_{i\pm1}e_i=e_{i},\ \ e_ie_j=e_je_i,\ {\rm
 if}\ |i-j|\ge 2,\cr
 &&t_i-t_{i}^{-1}=\epsilon(1-e_i),\cr
 &&e_it_i=t_ie_i=a e_i,\cr
 &&e_i^2=\tau e_i,\cr
 &&e_it_{i\pm 1}e_i=a^{- 1}e_i,\cr
 &&t_{i\pm1}t_ie_{i\pm1}=e_it_{i\pm 1}t_i=e_ie_{i\pm 1},
 \label{birman-wenzl}
\end{eqnarray}
with $a=p^{-1}q^{-2}$ and:
\begin{eqnarray}
 a&=&p^{-1}q^{-2},\ \ \epsilon=p-p^{-1},\cr
 \tau&=&1-{a-a^{-1}\over
 \epsilon} ={pq- p^{-1}q^{-1}\over p-p^{-1}}{(q+q^{-1})}.
 \label{definition-tau}
\end{eqnarray}
Notice that there is an algebra isomorphism from
$\mathcal{B}_n(\tau,a,\epsilon)$ to
$\mathcal{B}_n(\tau,a^{-1},-\epsilon)$ defined by $e_i\to e_i,\
t_i\to t_i^{-1}$.

This algebra can be endowed with a trace defined by:
\begin{eqnarray}
 &&{\rm tr}(1)=1,\cr
 &&{\rm tr} (xt_n)={1\over \tau a} {\rm tr} (x) \ \ {\rm for} \ \ x\in
 \mathcal{B}_{n-1}, \cr
 &&{\rm tr} (xe_n)= {1\over \tau}{\rm tr} (x) \ \ {\rm for} \ \ x\in
 \mathcal{B}_{n-1}.
 \label{trace-birman}
\end{eqnarray}

The B.M.W. algebra $\mathcal{B}_n$ has a natural representation in
terms of tangles, modulo Kaufman skein relations
\cite{kaufman}\cite{goodman}. Consider a rectangle with $2n$
points at its boundary. The $n$ points denoted $1,2,...,n$ are
disposed from left to right on the lower side and the $n$ points
denoted $\bar 1,\bar 2,...,\bar n$ are disposed above $1,2,...,n$
on the upper side. A tangle is made of $n$ strands projected onto
the rectangle connecting pairwise the $2n$ points at the boundary.
We denote the strand connecting $k$ to $l$ by ${k\choose l}$.

The tangles $\pi_1,\pi_2$ are composed by placing the rectangle
$\pi_2$ on top of $\pi_1$ so as to obtain a rectangle where the
lower points $k$ of $\pi_2$ coincide with the upper points $\bar
k$ of $\pi_1.$ The tangle denoted $\pi_1*\pi_2$ is obtained by
joining into a single stand the strand  of $\pi_1$ ending at $\bar
k$ with the strand of $\pi_2$ ending at $ k.$

The identity tangle is made of $n$ vertical strands connecting
$\bar k$ to $k$ :

$$Id=\pmatrix{\bar 1&\bar 2&...&\bar n\cr 1&2&...&n}$$.

The generator $t_k$ is obtained from the identity by exchanging
the two ends $k$ and $k+1$ in such a way that the strand ${\bar
k\choose k+1}$ crosses over ${\overline{k+1}\choose k}$.
Similarly, the generators $e_k$ is made of $n-2$ strands
connecting $l$ to $\bar l$ for $l\ne k,k+1$, two disjoint strands,
one above connecting $\bar k$ to $\overline {k+1}$, and one below
connecting $k$ to $k+1$.

The defining relations of the algebra (\ref{birman-wenzl}) can be
recast into skein relations allowing to identify tangles. A tangle
is multiplied by a factor $a^{\pm1}$ each time a loop is untwisted
and a close loop can be removed by multiplying the tangle by a
factor $\tau$. The relation $t_i-t_i^{-1}=\epsilon(1-e_i)$ enables
to relate over-crossings to under-crossings. The trace closes the
rectangle into a cylinder by identifying the upper with the lower
edge.

\subsection{Representation of the B.W.M. algebra on words}

We consider the case $n$ even. A Hilbert space $\mathcal{H}_n$ is
defined by acting with  $\mathcal{B}_{n}$ to the left on the minimal
projector $\alpha=e_1e_3...e_{n-1}$.

A basis is given by the reduced words:
 \begin{eqnarray}
 \pi=(t^{-1}_{c_{n-1}}t^{-1}_{c_{n-1}+1}...e_{n-1})...(
 t^{-1}_{c_{2p+1}}t^{-1}_{c_{2p+1}+1}...e_{2p+1})...e_1,
 \label{basis}
\end{eqnarray}
with $0\le p\le {n\over 2}-1$ and $1\le c_{2p+1}\le 2p+1$, and
each sequence reduces to $e_{2p+1}$ when $c_{2p+1}=2p+1$. There is
an order relation: $\pi>\pi'$ if $\pi$ is written $b\pi'$ with $b$
a monomial in $t_i^{-1}$.

Basis states can be  identified  with their tangle. The upper part
of the tangle is made of ${n\over 2}$ arches connecting
$\overline{2k-1}$ to $\overline{2k}$ and can be ignored. The lower
$n$ points $1,2,...,n$ are connected pairwise in all possible
ways. When two strands, ${a_1\choose b_1}$ and ${a_2\choose b_2}$
with $a_1<b_1,$ $a_2<b_2$ and $ b_1<b_2$ cross, ${a_2\choose b_2}$
crosses over ${a_1\choose b_1}$. A basis state is therefore
characterized by the sequence of strands:

 \begin{eqnarray}
\pi = \pmatrix{a_1&a_2&...&a_n\cr b_1&b_2&...&b_n},
 \label{sequence-basis}
\end{eqnarray}
with $a_1<a_2<...<a_n$, $a_i<b_i$.

The basis states (\ref{basis}) can be recursively identified with
tangles as follow: The tangle representing the word
$\pi=(t^{-1}_{c_{n-1}}...e_{n-1})\pi'$ with $\pi'\in
\mathcal{H}_{n-2}$ is obtained by pushing by one unit to the right
the strand ends of $\pi'$ larger or equal to $c_{n-1}$, and by
inserting a strand ${c_{n-1}\choose n}$ upon the tangle obtained
in this way.

Consider the maximally crossed state $\rho$ characterized by the
condition $t_i^{-1}\rho=t_{{n\over 2}+i}^{-1}\rho$ for $1\le i\le
{n\over 2}-1$ :
 \begin{eqnarray}
\rho = \pmatrix{1&2&...&{n\over 2}\cr {n\over 2}+1&{n\over
2}+2&...&n}.
 \label{maxicrossed-basis}
\end{eqnarray}
Another way to obtain the basis states is to act on $\rho$ with
the generators $t_i^{\pm 1},$ imposing that the number of
crossings decreases by one unit each time one acts with $t_i^{\pm
1}.$ The possible actions are given by:

 \begin{eqnarray}
 t_i^{-1}\pmatrix{i&i+1\cr j&k}&=&\pmatrix{i&i+1\cr k&j} \ {\rm if}\
 k>j,\cr
 t_i\pmatrix{j&...&i\cr i+1&...&k}&=&\pmatrix{j&...&i+1\cr i&...&k}\cr
 t_i^{-1}\pmatrix{j&...&k\cr i&...&i+1}&=&\pmatrix{j&...&k\cr i+1&...&i}.
\label{maxicrossed-action}
\end{eqnarray}

A different order relation than the preceding one is now
$\pi>\pi'$ if $\pi$ is written $b\pi'$ with $b$ a monomial in
$t_i^{\pm 1}.$ Thus one must be able to go from $\pi'$ to $\pi$ by
a sequence of untwists. The only minimal tangle is $\rho$ and all
the basis states can be reached from it by this procedure. We call
this order the crossing order to differentiate it from the first
one.

\smallskip

By identifying the two vertical sides of the rectangle, we can
view the labels $1,2,...,n$ as arranged cyclically around a
circle. We define $\sigma$, the operator which acts by cyclically
permuting the indices of a tangle $i\to i-1$:

\begin{eqnarray}
\sigma=at_n^{-1}t_{n-1}^{-1}...t_1^{-1}.
 \label{sigma}
\end{eqnarray}
The normalization factor can be determined by a straightforward
check that $\sigma^{-2}\alpha=\alpha$.

Note that  the application of $\sigma$ on a reduced word $\pi$
spoils the crossing rules and the relations (\ref{birman-wenzl})
need to be used to put the word in its canonical form back.
Nevertheless, $\sigma$ is a triangular matrix for the crossing
order: $\sigma \pi=  \pi_\sigma +{\rm \ tangles\ >}\pi$, where
$\pi_\sigma$ is the basis state obtained from $\pi$ by the
substitution $i\to i-1$. It follows from this that the line $\rho$
of the matrix representing $\sigma$ has only zeros except a $1$
for the diagonal element.

A Hermitian conjugation  is defined by $t_i^*=t_i^{-1}, \
e_i^*=e_i,\ p^*=p^{-1},q^*=q^{-1},$ and a scalar product as:
\begin{eqnarray}
\pi^*\pi'=\langle \pi|\pi'\rangle\alpha.
 \label{scalar}
\end{eqnarray}
\subsubsection{Structure of the matrices $e_i,t_i$ and dual basis}

 In the basis (\ref{basis}), the matrices $e_i,t_i$ have the
following structure. The basis states can be split into the
tangles $\Pi^i_0$ with a strand ${i\choose i+1} $ connecting $i$
to $i+1$, the tangles $\Pi^i_{1}$ where the two strands ending at
$i$ and $i+1$ cross and the tangles $\Pi^i_{2}$ where they do not
cross: $ t_i\Pi^i_{2}=\Pi^i_{1}$. We identify the tangles with the
projector onto these tangles,
$\Pi^i_0\oplus\Pi^i_1\oplus\Pi^i_2=1$. Thus,
$(\tau-e_i)\Pi^i_0=(t_i-a)\Pi^i_0=0$ and $(\Pi^i_{1}+\Pi^i_{2})
e_i=0.$ Using the relation $e_i(t_i-a)=0,$ the matrices $t_i,e_i$
are given the following block matrix expression:
\begin{eqnarray}
 e_i=\pmatrix{ \tau&a v&v\cr 0&0&0\cr 0&0&0},\ \
 t_i=\pmatrix{ a&-a\epsilon v&0\cr 0&\epsilon&1\cr 0&1&0},
 \label{matriceti}
\end{eqnarray}
with,
\begin{eqnarray}
 v=\Pi^i_0e_i\Pi^i_{2}.
 \label{matricev}
\end{eqnarray}

Let us now consider the dual basis of (\ref{basis}). Its basis
elements $\bar\pi$ are labeled by reduced words $\pi$ and are
defined by the relation:
\begin{eqnarray}
\bar\pi(\pi')=\delta_{\pi,\pi'}.
 \label{duality1}
\end{eqnarray}

If $\psi$ is an element of the B.M.W algebra, we have the duality
relation:
\begin{eqnarray}
\bar\pi\psi(\pi')=\bar\pi(\psi\pi'),
 \label{duality2}
\end{eqnarray}
so that we can generate the whole dual basis upon acting with the
B.M.W. algebra on the generator $\bar \alpha$. The matrices $t_i$
(\ref{matriceti}) are such that the column indexed by the tangle
$\pi\in\Pi^i_2$ has only zeros except for a $1$ at the
intersection with the line indexed by $t_i\pi\in\Pi^i_{1}$.
Similarly, the line $t_i\pi\in \Pi^i_{1}$ has only zeros except
for a $1$ at the intersection with the column $\pi\in\Pi_2^i$.
>From this, using the fact that the the basis elements are obtained
by successive action of the generators $t_i$ on $\Pi_2^{i}$, we
deduce that $\bar\pi$ can be identified with the reduced
expression of $\pi$ written in the reverse order:
\begin{eqnarray}
 \bar\pi=\bar\alpha...(t^{-1}_{2p}...t^{-1}_{a_{2p+1}+1}t^{-1}_{a_{2p+1}})...(
 t_{n-2}^{-1}...t_{a_n}^{-1}).
 \label{expressionFpi}
\end{eqnarray}
A similar argument shows that the same holds for words constructed
using $\bar{\rho}$ instead of $\bar{\alpha}$ as a generator.

\section{Polynomial representations}

\subsection{Projectors from the Yang-Baxter equation \label{yang-baxter-section}}

Let us consider the Yang-Baxter equation:
\begin{eqnarray}
 Y_{12}(z)Y_{23}(z w)Y_{12}(w)=Y_{23}(w)Y_{12}(z w)Y_{23}(z).
\label{yang-baxter0} \end{eqnarray}

It is possible to define an operator $Y_{ii+1}(z)$ solution of the
Yang-Baxter equation in terms of the algebra (\ref{birman-wenzl})
\cite{Jones2}\cite{ge}:
\begin{eqnarray}
 D(z)Y_{ii+1}(z)=(z-1)p^{-1}a^{-1}t_i-\epsilon(1-p^{-1}a^{-1})+(z^{-1}-1)t_i^{-1},
 \label{Y-birman-def}
 \end{eqnarray}
where the normalization factor $D(z)$ is chosen so that
$Y({z^{-1}})Y(z)=1.$ A possible choice for $D(z)$ is to require that
$Y_{ii+1}e_i=1$. Then, $D(z)$ is factorized as
$D(z)=(z-q^2)(z-p^2)p^{-1}z^{-1}$.

This allows us to identify consistently the Yang-Baxter matrix
$Y_{ii+1}$ with the permutation of the variables $z_i,z_{i+1}$ of a
vector $\Psi(z_1,...,z_n)$ parameterized by $z_i$:

\begin{eqnarray}
Y_{ii+1}({z_{i+1}\over z_{i}})\Psi=\Psi k_{ii+1}.
 \label{Y=k}
\end{eqnarray}
where  $k_{ii+1}$ acts by permuting the variables $z_i$ and
$z_{i+1}$ of the expressions located left of it.

We want to use the relation (\ref{Y=k}) to constrain operators
$\bar t_i$ and $\bar e_i$ which act on the components $\bar \pi$
of $\Psi$ and are dual to the $t_i,e_i$:
\begin{eqnarray}
t_i\Psi=\Psi \bar t_i,\ \ e_i\Psi=\Psi \bar e_i.
 \label{bar=pas-bar}
\end{eqnarray}

(\ref{Y=k}) does not allow to obtain a closed expression for $\bar
t_i$ and $\bar e_i$, but gives only a partial information about
them. By multiplying (\ref{Y=k}) by $e_i$ to the left, we obtain:
$e_i\Psi(k_{ii+1}-1)=0,$ which, with (\ref{bar=pas-bar}), implies
that $\bar e_i$ projects onto polynomials symmetric under the
exchange of $z_i,z_{i+1}$.

If we multiply (\ref{Y=k}) by $e_i-\tau$ to the left, we obtain:

\begin{eqnarray}
 (e_i-\tau)\Psi k_{ii+1}={qz_{i+1}-q^{-1}z_{i}\over
qz_{i}-q^{-1}z_{i+1}}{(z_{i+1}-z_{i})t_i+(p-p^{-1})z_i\over
pz_{i}-p^{-1}z_{i+1}}(e_i-\tau)\Psi.
 \label{e-tau(k-1)=0}
\end{eqnarray}
Since the left hand side is a polynomial, this implies that
$e_i-\tau$ projects $\Psi$ onto its components divisible by
$qz_{i}-q^{-1}z_{i+1}$. It is useful to introduce the notation
$b_{ij}=qz_i-q^{-1}z_{j}$. Thus,
\begin{eqnarray}
 \bar\psi\bar e_i=\tau\bar\psi&\Leftrightarrow&\bar\psi\
 {\rm is\ symmetrical\ under\ the\ exchange\ of}\ z_i,z_{i+1}\cr
 \bar\psi \bar e_i=0 &\Leftrightarrow& \bar\psi \ {\rm is\ divisible\ by\ } b_{ii+1}.
 \label{resume-e-psi}
\end{eqnarray}

$\bar e_i$ projects  onto the symmetrical polynomials under the
exchange of $z_i$ and $z_{i+1}$, orthogonally to the polynomials
divisible by $b_{ii+1}$. If we restrict to the components
$\bar\psi \bar e_i=0$, the action of the projectors $\bar t_i-p$
and $\bar t_i+p^{-1}$ acquire the expressions:

\begin{eqnarray}
\bar t_i-p&=&-b_{ii+1}^{-1}{pz_{i+1}-p^{-1}z_{i}\over
z_{i+1}-z_{i}}(1+k_{ii+1})b_{ii+1},\cr \bar
t_i+p^{-1}&=&b_{ii+1}^{-1}(1-k_{ii+1}){pz_{i}-p^{-1}z_{i+1}\over
z_i-z_{i+1}}b_{ii+1}, \label{polynomeTL-gi}
\end{eqnarray}

Another way to view the action of the generators $t_i$ on
polynomials is to use the affine-Hecke relations:
\begin{eqnarray}
 z_i\bar t_j&=&\bar t_jz_i\ \ {\rm if}\ i\ne j,j+1 \cr
 z_{i+1}\bar t_i&=&\bar t_i^{-1}z_i
 \label{affine-hecke}
\end{eqnarray}
to commute $\bar t_i^{\pm 1}$ through the polynomial and to make
use of the relation $b_{ii+1}\bar t_i=-p^{-1}b_{ii+1}$ when
$b_{ii+1}$ is the only monomial left.

\subsection{Statement of the problem}

Let us consider a vector $\Psi$:
\begin{eqnarray}
\Psi=\sum_{\pi}\pi \bar{\pi}(z_i), \label{PDF}
\end{eqnarray}
constructed in the following way. The vectors $\pi$ are the basis
vector of $\mathcal{H}_n$ indexed by $n$ even on which the  B.M.W.
algebra acts to the left. $\bar{\pi}$ are certain basis elements
of a space $\mathcal{M}_n$ of homogeneous polynomials of degree
$n({n\over 2}-1)$ in the variables $z_1,z_2,...,z_n$. We want to
determine an action of the affine B.M.W. algebra on the
polynomials, generated by the operators $\bar t_i$ and $\bar
\sigma $ acting to the right and determine the coefficients
$\bar\pi$ in such a way that both actions give the same result on
the vector $\Psi$:
\begin{eqnarray}
\Psi\bar t_{i}&=& t_{i}\Psi \cr\Psi \bar\sigma&=&
\sigma\Psi,\label{duality}
\end{eqnarray}

Said differently, we look for a dual action of the affine B.M.W.
algebra acting on polynomials. Unless we specify it, we address
this problem for a generic value of the deformation parameters,
for which the B.M.W. algebra is semisimple.

We  verify that this problem has a solution in the case $n=4$ in the
appendix \ref{construction M4}.

\subsection{The property defining the polynomial space\label{poly}}

 We consider a space of polynomials in an infinite number of variables $z_i,$ where
$z_{i+n}$ is identified with $a^{-2}z_i,$ for $a=p^{-1}q^{-2}$.
These polynomials are constructed as linear combinations of the
monomials $z^\mu=z_1^{\mu_1}...z_n^{\mu_n}$ with a total degree
$|\mu|$ fixed. The cyclic operator $\bar \sigma$ acts as:
\begin{eqnarray}
\bar\psi\bar\sigma(z_i)=a^{2|\mu|\over n}\bar\psi(z_{i+1}),
\label{sigma-pfaf}
\end{eqnarray}
where the normalization is such that $\bar \sigma^n=1$.

 We require that these polynomials obey the property (P):
\begin{eqnarray}
 {\rm(P): }\
&&\bar\psi(z_i=z,z_j=q^2z,z_k=q^4z)=0,\ \ {\rm if}\ i,j,k,\ {\rm
are\ cyclically\
 ordered}.
 \label{propiete-pfaf}
\end{eqnarray}
In order for (P) to be compatible with the cyclic identification,
$\bar\psi$ must vanish for the triplets:
$(z_i,z_j,z_k)=(z,q^2z,q^4z),\ (z,p^2z,q^2p^2z)\ {\rm\ and\ }
(z,q^2z,q^2p^2z)$ when $1\le i<j<k \le n$.

$\mathcal{M}_n$ is the space of minimal degree polynomials
satisfying (P). We shall show that a polynomial belonging to
$\mathcal{M}_n$ can be split in a unique way as:
\begin{eqnarray}
 \bar \psi&=&S_0+b_{ii+1}(S_1+(pz_i- p^{-1}z_{i+1})S_2),
\label{decompose-pfaf}
\end{eqnarray}
where the three polynomials in the decomposition belong to
$\mathcal{M}_n$, and $S_0,S_1,S_2$ are symmetrical under the
exchange of $z_i$ and $z_{i+1}$. To understand the origin of the
second projector, assume that $\bar \psi$ is divisible by
$b_{ii+1}$: $\bar \psi=b_{ii+1}\bar \phi$ . It is then possible to
act on $\bar\psi$ with the projectors (\ref{polynomeTL-gi}),
preserving the property (P). The factor $(qz_i-
q^{-1}z_{i+1})(pz_i- p^{-1}z_{i+1})$ in the definition of
$\bar\psi(\bar t_i+p^{-1})$ ensures that this polynomial vanishes
when $z_{i+1}=qz_i,$ or $z_{i+1}=pz_i,$ and therefore obeys (P).

The two projectors (\ref{polynomeTL-gi}) enable to achieve the
decomposition of the factor proportional to $b_{ii+1}$ in
(\ref{decompose-pfaf}). Thus, at this point we can write:
\begin{eqnarray}
 \mathcal{M}_n&=&\mathcal{M}_n^0\oplus\mathcal{M}_n^1\oplus\mathcal{M}_n^2\cr
 &=&S_0+b_{12}(S_1+(pz_1-
 p^{-1}z_{2})S_2),
\label{decompose-pfaf1}
\end{eqnarray}
with $S_1,S_2$ symmetrical under the exchange of $z_1,z_2,$ and $
\mathcal{M}_n^0$ is a complementary space to the polynomials in
$\mathcal{M}_n$ divisible by $b_{12}=qz_1-q^{-1}z_2$ . Since
$S_0\in \mathcal{M}_n^0$ is defined up to a polynomial divisible
by $b_{12}$, it remains to prove that it can be chosen symmetrical
under the exchange of $z_1,z_2$. The  section \ref{bare1} is
devoted to prove this fact.

By applying (P) to $z_1,z_2,z_i$, we define a projection from
$\mathcal{M}_n\to \mathcal{M}_{n-2}$ as:
\begin{eqnarray}
E'(\bar \psi)=c'{1\over \phi(z,z_i)}\bar
\psi(z_1={z},z_2=q^2z,z_i),\label{definitE'}
\end{eqnarray}
where  $\phi(z,z_i)=\prod_{i=3}^{n}(q^4z-z_i)(p^2q^2z-z_i),$ and
$c'$ is a normalization constant.
$E'(\mathcal M_{n})\subset\mathcal M_{n-2}$, and using a recursion
argument on $n$, the degree of the polynomials in $\mathcal M_{n}$
is $n({n\over 2}-1)$.

\bigskip

Let us construct a basis of polynomials in correspondence with the
dual states $\bar{\pi}$ (\ref{duality1}). A state obeying (P) with
the correct degree is the maximally crossed state given by:
\begin{eqnarray}
 \bar {\rho}=\prod_i \prod_{k=1}^{{n\over 2}-1} b_{ii+k},
 \label{crossed-pfaf}
\end{eqnarray}
with the cyclic identifications understood and
$b_{ij}=qz_i-q^{-1}z_j$.

One has $\bar\rho t_i=\bar \rho t_{{n\over 2}+i}$, it is therefore
consistent to construct a basis of polynomials $\bar\pi$, by
acting with the operators $\bar t_i,\bar t_i^{-1}$ on $\bar \rho$
as in (\ref{maxicrossed-action}).

Since the generators $\bar t_i,\bar t_i^{-1}$
(\ref{polynomeTL-gi}) act on polynomials divisible by $b_{ii+1}$,
we must verify that the polynomial upon which one acts with $\bar
t_i^{\pm 1}$ is always divisible by $b_{ii+1}$. We show here a
more general result.

If a tangle $\pi$ has consecutive points $j,j+1,...,k$ with no
strands connecting any two among them, we say that $j$ and $k$ are
matched, if not they are split. In particular, two consecutive
points are either matched, or connected by a strand. Let us show
that when $j$ and $k$ are matched, the polynomial $\bar \pi$ is
divisible by $b_{jk}.$

The property is true for $\bar\rho$, and we can use a recursion
argument on the crossings. If two points are split in $\pi<\pi'$,
they are split in $\pi',$ otherwise the number of crossings would
not decrease continuously when going from $\pi$ to $\pi'$. The
only way $t_i$ splits $j<k$ from $k$ is when $i=j-1,$ or $i=k.$
Using the fact that $t_i$ commutes with $b_{ik}b_{i+1k}$
($b_{ki}b_{ki+1}$) as well as $b_{kl}$ for $k$ and $l\ne i,i+1$,
we deduce that when $i$ is matched with $k$ the polynomial
$\bar\pi$ is divisible by $b_{ik}$.

Once we have identified a basis of polynomial with the dual basis
of $\mathcal{H}_n$, we have determined an action of
$\mathcal{B}_{n}$ on $\mathcal{M}_n$. We still need to verify that
this action is consistent with the decomposition
(\ref{decompose-pfaf}), and in particular that $e_i\Psi$ is
symmetrical under the exchange of $z_i$ and $z_{i+1}.$

The basis elements can be ordered by their highest monomial
$z^{\lambda}$ (see \cite{incompressible}), and this ordering is
compatible with the first ordering defined on the tangles. The
highest monomial of the basis is given by
$\lambda=(n-2,n-2,n-4,n-4,...,0,0).$

\subsection{The projectors $\bar e_i$ and the full duality relation \label{bare1}}

We obtain an expression of $\bar e_1$ as an operator acting on the
components $\Pi_0^1$ having the two indices $1,2$ connected by a
link (split). The argument is essentially an adaptation to this
case of the proof of Di Francesco and Zinn-Justin \cite{Pdf2}.

>From the definition (\ref{resume-e-psi}) of $\bar e_1$, we know
that it is null on $\Pi^1_{1,2}$. Consider $\bar \pi_0\in
\Pi^1_0$. The action of $\bar e_1$ on $\bar \pi_0$ is fully
determined by the condition that $\bar \pi_0 \bar e_1=\tau \bar
\pi_0 +\bar \pi_1$ is symmetrical under the exchange of $z_1,z_2,$
and $\bar \pi_1$ is divisible by $b_{12}.$

We define the  projector:
\begin{eqnarray}
 \bar f_1=(1-k_{12}){1\over
z_{1}-z_{2}}(pz_1-p^{-1}z_2)(qz_1-q^{-1}z_2){q+q^{-1}\over
p-p^{-1}}. \label{polynomeTL-ei-wenzl}
\end{eqnarray}
$\bar f_1$ raises the degree of the polynomial upon which it acts
by one, its kernel are the symmetric polynomials in $z_1,z_2$, and
$\bar f_1^2= (z_1+z_2)\tau \bar f_1$. The complementary projector:
$(z_1+z_2)\tau-\bar f_1$ projects onto a polynomial symmetrical
under the exchange of $z_1,z_2$.

Consider the following tangle $\pi_0 \in \Pi_0^1:$

 \begin{eqnarray}
 \pi_0=\pmatrix{1&3&...&{n\over2}+1\cr
 2&{n\over2}+2&...&n}&=& t_2t_3...t_{n\over2}\rho\cr
 =ae_1t_2^{-1}...t_{i-1}^{-1}t_{i+1}...t_{n\over2}\rho
 &=&e_1t_1^{-1}t_2^{-1}...t_{i-1}^{-1}t_{i+1}...t_{n\over2}\rho,
 \label{exemple-1-birman}
\end{eqnarray}
for $2\le i\le {n\over 2}.$ With the diagonal element
$e_{\pi_0\pi_0}=\tau,$ the second line of (\ref{exemple-1-birman})
exhausts all the matrix elements $e_{\pi_0\pi}$ of the line
$\pi_0$ of $e_1$. By duality, we must have:

 \begin{eqnarray}
 \bar\pi_0\bar e_1=\sum_\pi e_{\pi_0\pi}\bar \pi.
 \label{exemple-1-birman-dual}
\end{eqnarray}
If (\ref{resume-e-psi}) holds, the right hand side of this
equality must be a symmetrical polynomial in $z_1,z_2$. If we set:
\begin{eqnarray}
\bar\pi_2=\sum_{i=2}^{n\over 2}\bar\rho \bar t_{n\over 2}...\bar
t_{i+1}\bar t_{i-1}^{-1}...\bar t_2^{-1},
 \label{exemple-11-birman}
\end{eqnarray}
this amounts to show that:
\begin{eqnarray}
 \tau \bar \pi_0+\bar \pi_2(a+\bar t_1^{-1}),
 \label{dual-example-2-birman}
\end{eqnarray}
is a symmetrical polynomial in $z_1,z_2$.

Let us consider the polynomial:
 \begin{eqnarray}
 \bar\pi'=\bar\pi_0 \bar f_1+\bar \pi_2(a+\bar t_1^{-1})(z_1+z_2)+q^{-1}\bar\pi_2(p-\bar
 t_1)b_{21}.
 \label{dual-example-3-birman}
\end{eqnarray}
The third term of the above sum is symmetrical from the definition
(\ref{polynomeTL-gi}) of $p-\bar t_1$. From the fact that $\bar
f_1-(z_1+z_2)\tau$ projects on a symmetric polynomial, the sum of
the two first terms, and thus the hole expression, is symmetrical
if and only if (\ref{dual-example-2-birman}) holds. After
substituting the expressions (\ref{polynomeTL-gi}) of $\bar t_1$
and (\ref{polynomeTL-ei-wenzl}) of $\bar f_1$, the sum
(\ref{dual-example-3-birman}) can be put under the form
$\bar\pi'=S_n\bar f_1$ with:
 \begin{eqnarray}
 S_n=\bar \pi_0+\epsilon q^{-1} b_{12}^{-1}z_2\bar \pi_2.
 \label{dual-example-4-birman}
\end{eqnarray}
>From the definition of $\bar f_1$, the only way  $S_n\bar f_1$ can
be symmetrical is when it is equal to zero. Since the kernel of
$\bar f_1$ are the symmetric polynomials, we therefore need to
show that $S_n$ (\ref{dual-example-4-birman}) is symmetrical under
the exchange of $z_1$ and $z_2$. Denoting $b_{j}=b_{1j}$, and
substituting the explicit expressions (\ref{exemple-1-birman},
\ref{exemple-11-birman}) of $\bar \pi_0$, $\bar \pi_2,$ we obtain:

 \begin{eqnarray}
 S_n=\bar\rho (\bar t_{n\over2}...\bar t_3\bar t_2b_2+
 \epsilon q^{-1}\sum_{i=2}^{n\over 2} \bar t_{n\over 2}...\bar
 t_{i+1}\bar t_{i-1}^{-1}...\bar t_2^{-1}z_2) b_2^{-1}=\bar\rho
 I_n.
 \label{dual-expression-birman}
\end{eqnarray}

A first application of the affine Hecke relations
(\ref{affine-hecke}) enables to replace $\bar t_{i-1}^{-1}...\bar
t_2^{-1}z_2$ with $z_i\bar t_{i-1}...\bar t_2$ in
(\ref{dual-expression-birman}). A second one gives:
\begin{eqnarray}
 \bar t_i b_i+\epsilon q^{-1}z_i=b_{i+1}\bar t_i,
 \label{dual-intermediaire-birman}
\end{eqnarray}
for $i\ge 2.$ The repeated application of this identity enables to
put the operator $I_n$ acting on $\bar \rho$ in
(\ref{dual-expression-birman}) under the form:
  \begin{eqnarray}
I_n= b_{{n\over 2}+1}\ \bar t_{n\over 2}...\bar t_{2}\
 b_{2}^{-1},
 \label{dual-final-birman}
\end{eqnarray}
Using the fact that $\bar t_i$ commutes with $b_{ij}b_{i+1j},$ and
its action on the product $b_{i+1j}b_{ii+1}$:
 \begin{eqnarray}
 b_{i+1j}b_{ii+1}\bar
 t_i=b_{ji+n}b_{ii+1},
 \label{action-rho-brauer}
\end{eqnarray}
we obtain a factorized expression of $S_n$
(\ref{dual-expression-birman}) exhibiting the symmetry under the
exchange of $z_1,z_2$:
 \begin{eqnarray}
S_n=\bar\rho_{n-2}(\prod_{k=3}^{{n\over
2}+1}b_{1k}b_{2k})(\prod_{k={n\over 2}+2}^n b_{kn+1}b_{kn+2}),
 \label{action-factorisee-birman}
\end{eqnarray}
where $\bar \rho_{n-2}$ is the maximally crossed state of
$\mathcal{M}_{n-2}$ expressed in the variables $z_i,\ 3\le i\le
n$.

Thus, the relation (\ref{exemple-1-birman-dual}) holds. Let
$\mathcal{B}_{n-2}$ be the subalgebra generated by $t_i$, $i\ge
3.$ Once we know that the line $\pi_0$, of the matrices $\bar e_1$
and $e_1$ are the same, acting with $\mathcal{B}_{n-2}$ on
$\bar\pi_0$ and using the property that $\mathcal{B}_{n-2}$
commutes with $e_1$, we can equate the lines  $\pi\in\Pi_0^1$.
This establishes that $e_1$ projects $\Psi$ onto symmetrical
polynomials under the exchange of $z_1$ and $z_2$ orthogonally to
the polynomials divisible by $b_2$.

>From the matrix expression (\ref{matriceti}), we can extend the
definition of $\bar t_1=t_1$ to $\Pi_0^1$, and thus to the
polynomials of $\mathcal{M}_n$ not divisible by $b_2$. We  define
$\bar t_i$ for $i\ge 1$ through $\bar \sigma\bar t_i=\bar
t_{i-1}\bar\sigma ,$ and this definition agrees with
(\ref{polynomeTL-gi}) on polynomials divisible by $b_{ii+1}$.
$\bar t_i$ and $\bar e_i$ constructed with this procedure are
identified with $t_i$ and $e_i$

To complete the identification we must verify that $\bar\sigma$
(\ref{sigma-pfaf}) is equal to $\sigma$ (\ref{sigma}). We observe
that $\bar\rho\bar\sigma=\bar\rho$, in agreement with the fact
that the line $\rho$ of the matrix $\sigma$ has only zeros except
a $1$ for the diagonal element. Since $y_1=\sigma^{-1}\bar \sigma$
commutes with $\mathcal{B}_{n-1}$, the algebra generated by $t_i$
for $i\ge 2$, and $\mathcal{H}_n$ is an irreducible module over
$\mathcal{B}_{n-1}$, we deduce that $y_1=1$.

\subsubsection{The cyclic case}

In the cyclic case $p^{-1}=q^2$, following the line of
\cite{Pdf2}, it is straightforward to verify that $\Psi$ is the
perron-Frobenius eigenvector of a statistical transfer matrix. The
Pfaffian state (\ref{Pfaffien-pfaf}) is a symmetrical polynomial
obeying the constraint (\ref{propiete-pfaf}), and therefore
realizes a trivial representation of the B.M.W algebra.

We can verify this property in the tangle representation. Another
way to write the basis elements (\ref{basis}) is by dressing the
Temperley-Lieb words with braid generators. The basis is be
obtained from the Temperley-Lieb basis \cite{incompressible} by
substituting the $e_i$ for $i<d_{2p+1}$ sitting at the end of the
sequence $e_{c_{2p+1}}...e_{2p+1}$ with $t^{-1}_i$:
 \begin{eqnarray}
 \pi=(t^{-1}_{c_{n-1}}t^{-1}_{c_{n-1}+1}...e_{d_{n-1}}...e_{n-2}e_{n-1})...(
 t^{-1}_{c_{2p+1}}t^{-1}_{c_{2p+1}+1}...e_{d_{2p+1}}...e_{2p}e_{2p+1})...e_1,
 \label{basis2}
\end{eqnarray}
where $c_{2p+1}\le d_{2p+1} \le 2p+1$ and
$c_{n-1}>c_{n-2}>...>c_3>1$. In this presentation, it is clear
that when $a=\tau=1$, the generators $e_i$ transforms a basis
element into another basis elements. Therefore, the matrix
elements of each columns of the matrix $v$ (\ref{matricev}) are
all equal to zero except for a single one which is equal to one.
Consequently, the sum of dual vectors $\sum_\pi \bar\pi$ form a
trivial representation and is proportional to
(\ref{Pfaffien-pfaf}).

\bigskip
\bigskip

\appendix

\section {Explicit construction of $\mathcal{M}_4$ \label{construction M4}}

Let us construct $\mathcal{M}_4$ the dual of $\mathcal{H}_4$.

We search for an invariant vector $\Psi$ of the form:

\begin{eqnarray}
\Psi=\bar\pi_1(z_1,..,z_4) e_1e_3+
\bar\pi_2(z_1,..,z_4)t_2^{-1}e_1e_3+\bar\pi_3(z_1,..,z_4)e_2e_1e_3.\label{birman-vecteur-H4}
\end{eqnarray}

The dual representation is obtained by acting with the generators
on the maximally crossed state $\bar\pi_2\equiv(0,1,0)$ :

\begin{eqnarray}
\bar\pi_2=(qz_1-q^{-1}z_2)(qz_2-q^{-1}z_3)(qz_3-q^{-1}z_4)(qz_4-p^2q^{3}z_1).\label{birman-F3-H4}
\end{eqnarray}
>From (\ref{polynomeTL-gi}), we obtain
$\bar\pi_3\equiv(0,0,1)=\bar\pi_2\bar
 t_1^{-1}$ and $\bar\pi_1\equiv (1,0,0)=\bar\pi_2\bar
 t_2$:
\begin{eqnarray}
 \bar\pi_3&=&p^{-1}(p^2z_3-z_4)(q^2z_3-z_4)(p^2z_1-z_2)(q^2z_1-z_2)-p\bar\pi_2,\cr
 \bar\pi_1&=&q^{-2}p^{-1}(q^4z_1-z_4)(p^2q^2z_1-z_4)(p^2z_2-z_3)(q^2z_2-z_3)-p^{-1}\bar\pi_2.
 \label{birman-F12-H4}
\end{eqnarray}
>From these expressions, one obtains the expression of $\bar t_i$
using (\ref{polynomeTL-gi}) and of  $\bar e_i$ by projecting
$\bar\pi_j$ onto a polynomial symmetrical under the exchange of
$z_i$ and $z_{i+1}$.

The matrices representing the algebra $\mathcal{B}_4$
(\ref{birman-wenzl}) one obtains this way coincide with the tangle
representation and are given by:
\begin{eqnarray}
 e_1=e_3=\pmatrix{\tau&a&1\cr 0&0&0\cr 0&0&0}&,& \
 e_2=e_4=\pmatrix{0&0&0\cr 0&0&0\cr 1&a^{-1}&\tau},\cr
 t_1=t_3=\pmatrix{a&-a\epsilon&0 \cr 0&\epsilon&1\cr 0&1&0}&,& \
 t_2=t_4=\pmatrix{\epsilon&1&0\cr  1&0&0\cr -\epsilon&0&a}.
 \label{birman-matrix-H4}
 \end{eqnarray}

If we define the operator $\sigma$ acting by cyclically permuting
the indices of a tangle:
\begin{eqnarray}
 \sigma= a t_3^{-1}t_2^{-1}t_1^{-1}=\pmatrix{0&\epsilon&1\cr 0&1&0\cr
 1&-\epsilon&0},
 \label{sigma-H4}
\end{eqnarray}
one can verify that its action on $\Psi$ coincides with that of
$\bar\sigma$ (\ref{sigma-pfaf}): $\sigma\Psi=\Psi\bar\sigma$.

In the case $a=1$, $\tau=1$, and $\bar\pi_{\Omega}=(1,1,1)$ is a
trivial representation of the B.M.W. algebra.

\section {The Brauer Algebra \label{brauer-appendice}}

The Brauer algebra is the rational limit of the B.M.W. algebra
obtained in the limit $z=e^{2\epsilon u}, q=e^{\epsilon b},
p=e^{\epsilon c}$, with $\epsilon \to 0$. Then, $\tau=2(1+{b\over
c})$ and the relations (\ref{birman-wenzl}) become:
\begin{eqnarray}
 &&t_it_{i+1}t_i=t_{i+1}t_it_{i+1},\ \ t_it_j=t_jt_i,\ {\rm
 if}\ |i-j|\ge 2,\cr
 &&e_ie_{i\pm1}e_i=e_{i},\ \ e_ie_j=e_je_i,\ {\rm
 if}\ |i-j|\ge 2,\cr
 &&t_i^2=1,\cr
 &&e_it_i=t_ie_i=e_i,\cr
 &&e_i^2=\tau e_i,\cr
 &&e_it_{i\pm 1}e_i=e_i,\cr
 &&t_{i\pm1}t_ie_{i\pm1}=e_it_{i\pm 1}t_i=e_ie_{i\pm 1}.
 \label{brauer}
\end{eqnarray}

In this limit, the $t_i$ are transpositions and therefore, the
tangles become link patterns where one ignores crossings.

The structure of the matrices $t_i,e_i$ (\ref{matriceti}) is
simplified. In the formulas (\ref{matriceti}) , one can set $a=1$,
$\epsilon=0,$ $\tau=2(1+{b\over c})$ and all the nonzero matrix
elements of $v$ are equal to $1$.

The property (P) is replaced by:
\begin{eqnarray}
 {\rm(P'): }\
\bar\psi(u_i=u,u_j=u+b,u_k=u+2b)=0,\ \ {\rm if}\ i,j,k,\ {\rm are\
cyclically\
 ordered}.
 \label{propiete-brauer}
\end{eqnarray}

The cyclic identification becomes: $u_{i+n}=u_{i}+c+2b$ and it
imposes that  $\bar\psi$ vanishes for
$(x_i,x_j,x_k)=(u,u+b,u+2b),\ (u,u+c,u+c+b)\ {\rm\ and\ }
(u,u+b,u+c+b)$ when $1\le i<j<k\le n$.

The operator $Y_{ii+1}$ (\ref{Y=k}) which permutes the coordinates
is given by:

\begin{eqnarray}
 Y_{ii+1}(u)={(b-u)(c-ut_{i})+cue_{i}\over (b+u)(c+u)}.
 \label{Y-brauer}
\end{eqnarray}
$Y$ is a solution of the rational Yang-Baxter equation:
\begin{eqnarray}
 Y_{12}(u)Y_{23}(u+v)Y_{12}(v)=Y_{23}(v)Y_{12}(u+v)Y_{23}(u),
\label{yang-baxter-brauer}
\end{eqnarray}
and it is normalized so that $Y_{12}(u)Y_{21}(-u)=1$.

$\bar e_i$ projects  onto the symmetrical polynomials under the
exchange of $u_i$ and $u_{i+1}$, orthogonally to the polynomials
divisible by $b_{ii+1}=u_i-u_{i+1}+b$. If we restrict to the
components $\bar\psi \bar e_i=0$ , the action of the projectors
$\bar t_i-1$ and $\bar t_i+1$ acquire the expressions:

\begin{eqnarray}
\bar\phi(\bar t_i-1)&=&-b_{ii+1}^{-1}{u_{i+1}-u_i+c\over
u_{i+1}-u_{i}}(1+k_{ii+1})b_{ii+1},\cr (\bar
t_i+1)&=&b_{ii+1}^{-1}(1-k_{ii+1}){u_{i}-u_{i+1}+c\over
u_i-u_{i+1}}b_{ii+1}.\label{polynomeTL-gi-brauer}
\end{eqnarray}

Equivalently, the $t_i$ can be commuted through the polynomial
using the degenerate affine relations:
\begin{eqnarray}
 &&u_i\bar t_j=\bar t_ju_i\ \ {\rm if}\ i\ne j,j+1 \cr
  &&\bar t_i u_i-u_{i+1}\bar t_i=c,
 \label{degenerate-affine}
\end{eqnarray}
and their action on the identity is: $b_{ii+1}\bar t_i=-b_{ii+1}.$

The polynomials are no longer homogeneous and the maximally crossed
state is given by:
\begin{eqnarray}
 \bar \psi=\prod_i \prod_{0< k<{n\over 2}} b_{ii+k},
\label{Fomegabar-brauer}
\end{eqnarray}
with the cyclic identifications understood. The highest degree
monomial can easily be obtained in the limit $b=c=0$ and is equal
to:
\begin{eqnarray}
(-)^{\epsilon}{\prod_{i<j}(u_i-u_j)\over \prod_{i\sim j} (u_i-u_j)},
\label{highest-brauer}
\end{eqnarray}
Where $i\sim j$ means that $i$ and $j$ are linked, and
$(-)^{\epsilon}$ is the parity the number of transpositions $t_i$ in
the reduced expression of the word.

The polynomials can be decomposed as:
\begin{eqnarray}
 \mathcal{M}_n&=&\mathcal{M}_n^0\oplus\mathcal{M}_n^1\oplus\mathcal{M}_n^2\cr
 &=&S_0+b_{12}(S_1+(u_1-u_2+c)S_2),
\label{decompose-pfaf1}
\end{eqnarray}
with $S_0,S_1,S_2$ symmetrical under the exchange of $u_1,u_2.$ Note
that $\mathcal{M}_n^{0,2}$ have a highest degree monomial even under
the exchange of $u_1,u_2$, while $\mathcal{M}_n^{1}$ have an odd
highest degree monomial.

\smallskip
We obtain an explicit expression of $\bar e_1$ as an operator
acting on the components $\Pi_0^1\Psi=\Psi_0$ having the two
indices $1,2$ connected by a link (split). We know that the action
of $\bar e_1$ on $\Psi_{1,2}$ is equal to zero since these
components are divisible by $b_{12}$. We also use the property
that $\left(e_1+(\tau-e_1){1\over
2}(1+t_1)\right)\Psi_0=\tau\Psi_0$ which simply results from the
fact that the highest degree polynomial of $\Psi_0$ is even under
the exchange of $u_1,u_2$ as can be seen explicitly from
(\ref{highest-brauer}). Thus
$\Psi_0\subset\mathcal{M}_n^0\oplus\mathcal{M}_n^2$.

Let us consider the two complementary projectors defined as:

\begin{eqnarray}
 \bar e_1&=&{1\over c}b_{ii+1}{(u_2-u_1+c)\over u_2-u_1}(1+k_{12}),\cr
 \tau-\bar e_1&=&{1\over c}(1-k_{12}){(u_1-u_2+c)\over u_1-u_2}b_{ii+1}.
\label{polynomeTL-ei-brauer}
\end{eqnarray}
The normalization is fixed so that $e_1^2=\tau e_1$. When
restricted to $\mathcal{M}_n^0\oplus\mathcal{M}_n^2,$ $\bar e_1$
projects onto $\mathcal{M}_n^0$, the symmetrical polynomials under
the exchange of $u_1,u_2,$ orthogonally to $\mathcal{M}_n^2$. Thus
the action of $\bar e_1$ on $\Psi_0$ projects out the
$\mathcal{M}_n^2$ component and keeps the component of $\Psi_0\in
\mathcal{M}_n^0$.

We want to identify $\bar e_1$ with $e_1$ using (\ref{bar=pas-bar}).
We first verify this identification on a specific line of the matrix
$e_1$. For this, consider the following tangle $\pi_0 \in \Pi_0^1:$

 \begin{eqnarray}
 \pi_0=\pmatrix{1&3&...&{n\over2}+1\cr
 2&{n\over2}+2&...&n}&=& t_2t_3...t_{n\over2}\rho\cr
 =e_1t_2...t_{i-1}t_{i+1}...t_{n\over2}\rho&=&e_1t_1t_2...t_{i-1}t_{i+1}...t_{n\over2}\rho.
 \label{exemple-1}
\end{eqnarray}
We also have $e_1\pi_0=\tau\pi_0,$ and this exhausts all the matrix
elements $e_{\pi_0\pi}$ of the line $\pi_0$ of $e_1$. By duality, we
must have $\bar\pi_0\bar e_1=\sum_\pi e_{\pi_0\pi}\bar \pi$:

 \begin{eqnarray}
\bar\rho \bar t_{n\over2}...\bar t_3\bar t_2(\bar \tau-\bar
e_1)+\sum_{i=2}^{n\over 2}\bar\rho \bar t_{n\over 2}...\bar
 t_{i+1}\bar t_{i-1}...\bar t_2(1+\bar t_1)=0.
 \label{dual-example-1}
\end{eqnarray}

 Denoting $b_{j}=b_{1j}$, and substituting the expressions
(\ref{polynomeTL-gi-brauer},\ref{polynomeTL-ei-brauer}) of $\bar
t_1+1$ and $\tau-\bar e_1$, the above equality is equivalent to
require that:

 \begin{eqnarray}
 S_n=\bar\rho (\bar t_{n\over2}...\bar t_3\bar t_2+\sum_{i=2}^{n\over 2} \bar t_{n\over 2}...\bar
 t_{i+1}\bar t_{i-1}...\bar t_2 cb_2^{-1})=\bar\rho I_n
 \label{dual-expression}
\end{eqnarray}
is symmetrical under the exchange of $u_1,u_2$.

>From the degenerate affine Hecke relations
(\ref{degenerate-affine}), we have:
\begin{eqnarray}
 \bar t_ib_i+c=b_{i+1}\bar t_i.
 \label{dual-intermediaire-brauer}
\end{eqnarray}
This identity enables us to put the operator $I_n$ acting on $\bar
\rho$ in (\ref{dual-expression}) under the form:
  \begin{eqnarray}
I_n= b_{{n\over 2}+1}\ \bar t_{n\over 2}...\bar t_{2}\
 b_{2}^{-1}
 \label{dual-final-brauer}
\end{eqnarray}
Using the fact that $\bar t_i$ commutes with $b_{ij}b_{i+1j}$ and
its action on the product $b_{i+1j}b_{ii+1}$:
 \begin{eqnarray}
 b_{i+1j}b_{ii+1}\bar
 t_i=b_{ji+n}b_{ii+1},
 \label{action-rho-brauer}
\end{eqnarray}
we obtain a factorized expression of $S_n$ (\ref{dual-expression})
exhibiting the symmetry under the exchange of $u_1,u_2$:
 \begin{eqnarray}
S_n=\bar\rho_{n-2}(\prod_{k=3}^{{n\over
2}+1}b_{1k}b_{2k})(\prod_{k={n\over 2}+2}^n b_{kn+1}b_{kn+2}),
 \label{action-factorisee-brauer}
\end{eqnarray}
where $\bar \rho_{n-2}$ is the maximally crossed state of
$\mathcal{M}_{n-2}$ expressed in the variables $u_i,\ 3\le i\le
n$. Thus, the relation (\ref{dual-example-1}) dual to
(\ref{exemple-1}) holds.

\section { Temperley Lieb  and $k=2$ cancellation conditions \label{T.L.appendice}}

This appendix is an addendum to \cite{incompressible}, where it
was observed that in the $k=2$ case, the Hecke generators $e_i$
did obey the more constrained Temperley and Lieb algebra
relations. We show that this is a simple consequence of the
minimal degree hypotheses in the $k=2$ case.

We consider the space $\mathcal{M}_n$ of polynomials obeying the
condition (P)(\ref{propiete-pfaf}) and of the minimal possible
degree. We assume that the generators $\bar e_1,$ projecting
$\bar\psi\in\mathcal{M}_n$ onto a symmetrical polynomial under the
exchange of $z_1$ and $z_{2},$ obey the Hecke relations:
 \begin{eqnarray}
\bar e_1\bar e_{2}\bar e_1-\bar e_1=\bar e_{2}\bar e_{1}\bar
e_{2}-\bar e_{2}.
 \label{action-rho-brauer}
\end{eqnarray}
We denote by $U$ the above projector and we show here that the
Temperley and Lieb condition: $U=0,$ is realized when $k=2$.

$U$ projects $\bar\psi$ onto a  polynomial symmetric under the
permutations of $z_1,z_{2},z_{3}$ obeying (P). This implies that
$\bar \psi U$ vanishes if $z_1=z,z_2=q^2z,z_3=q^{-2}z,$ and
therefore, $E'(\bar \psi U)$ (\ref{definitE'}) is divisible by
$z_3-q^{-2}z$. Since $E'(\bar \psi U)$ does not depend on $z$, it
is equal to zero. The restriction of $E'$ to $\mathcal{M}_n^0$
being injective, this implies that $\bar \psi U=0$.

\end{document}